\documentclass[reqno]{amsart} %12pt,
\usepackage{amssymb}\usepackage{epsf}

\theoremstyle{definition}

 \numberwithin{equation}{section}

\newcommand{\be}{\begin{equation}}
\newcommand{\ee}{\end{equation}}
\newcommand{\ba}{\begin{eqnarray}}
\newcommand{\ea}{\end{eqnarray}}

\newcommand{\lab}[1]{\label{#1}}
\newcommand{\N}{\mathbb N}

\newcommand{\C}{\mathbb C}
\newcommand{\R}{\mathbb R}
\newcommand{\Z}{\mathbb Z}

\newcommand{\T}{\mathbb T}

\newcommand{\ve}{\varepsilon}

\begin{document}

\title[Continuous biorthogonal functions]
{Continuous biorthogonality of the \\
 elliptic hypergeometric function}

 \author{V. P. Spiridonov}
\address{Bogoliubov Laboratory of Theoretical Physics, JINR,
 Dubna, Moscow reg. 141980, Russia;
{\em e-mail}: {\tt spiridon@theor.jinr.ru}}

\thanks{To be published in {\em Algebra i Analiz (St. Petersburg Math. J.)}}

\begin{abstract}
We construct a family of continuous biorthogonal functions
related to an elliptic analogue of the Gauss
hypergeometric function. The key tools
used for that are the elliptic beta integral and the
integral Bailey chain introduced earlier by the author.
Relations to the Sklyanin algebra and elliptic
analogues of the Faddeev modular double
are discussed in detail.

\end{abstract}

\maketitle

\tableofcontents

\section{Introduction}

Classical and quantum completely integrable systems serve as
a rich source of special functions. Complexity of these systems
correlates with the complexity of solutions of the Yang-Baxter
equation -- rational, trigonometric, or elliptic \cite{tf}.
In the theory of special functions this structural hierarchy is
reflected in the existence
of the plain and $q$-hypergeometric functions \cite{aar}
and their elliptic generalizations \cite{spi:thesis}.
One of the aims of the present paper is to clarify
certain questions of such a correspondence at the
elliptic level.

In \cite{spi:aa,spi:thesis}, the author introduced an elliptic
analogue of the Gauss hypergeometric function $V(t_1,\ldots,t_8)$
generalizing many special functions of hypergeometric type with
``classical" properties \cite{aar}. In particular, it was shown that this
function exhibits symmetry transformations tied to the
exceptional root system $E_7$ and satisfies the elliptic hypergeometric
equation. We start this work with showing in Sect. 2 that the
$V$-function satisfies a simple biorthogonality
relation corresponding to the continuous values of a spectral
parameter. We use for that the elliptic beta integral \cite{spi:umn}
and an integral analogue of the Bailey chains introduced in \cite{spi:bai}.
Our construction can be considered as an integral generalization
of Rosengren's approach to the elliptic $6j$-symbols \cite{ros:elementary}.

In Sect. 3, we investigate the generalized eigenvalue problem for
a finite difference operator of the second order $D(a,b,c,d;p;q)$.
Connection of such problems with the biorthogonality relations was
considered in \cite{zhe:gevp}. A problem analogous to ours was
investigated earlier in \cite{rai:abelian,ros:elementary} for the discrete
spectrum. Its general solution defines a basis in some space
of meromorphic functions. For a simple scalar product, we derive
the dual basis defined by solutions of a similar spectral problem
for conjugated operators. As a result, the overlap of the basis
vectors with their duals coincides with the $V$-function.
In Sect. 4, we derive a new contiguous relation for the elliptic
hypergeometric function and the biorthogonality condition
 associated with this relation.

It turns out that the $V$-function is directly connected with the
Sklyanin algebra -- one of the central objects in the quantum
inverse scattering method \cite{skl1,skl2}. This relation
is established in Sect. 5 on the basis of the observation due
to Rains \cite{rai:trans,ros:sklyanin} that the $D$-operator
can be represented as a linear combination of all four
generators of the Sklyanin algebra. In this way we obtain a
particular infinite dimensional module of this algebra
with the continuous values of the Casimir operators.

It is well known that the Sklyanin
algebra can be degenerated to the quantum algebra $U_q(sl_2)$,
which was intensively investigated during the last two decades.
In \cite{fad:mod}, Faddeev has introduced the
modular double -- the quantum algebra $U_q(sl_2)\otimes U_{{\tilde q}^{-1}}(sl_2)$,
where $\tilde q$ is a modular transform of $q$ (see also \cite{kls:unitary}).
Originally, it was aimed at dealing with the well defined logarithms of quantum
algebra generators, which can be traced to the demand of analytical
uniqueness of the algebra representation modules. Applying the latter requirement
to the results of Sect. 5, we come in Sect. 6 to two different elliptic
analogues of the modular double defined by direct products of
two copies of the Sklyanin algebra with the generators related
to each other by permutation of certain parameters.
An interesting fact is that some of the generators
of these two algebras anticommute with each other.

Finally, in Sect. 7 we discuss the Bethe ansatz for the $D$-operator
eigenfunctions and relations to the Calogero-Sutherland type models --
common topics in the theory of integrable systems. We demonstrate that
the standard eigenvalue problem for the $D$-operator
can be interpreted as a difference analogue of the Heun equation.
The latter equation is known to emerge in the one particle sector of
the Inozemtsev integrable model \cite{ino:lax}. Our finite difference
analogue of it turns out to be related to the model introduced by van Diejen
\cite{die:integrability} for a restricted set of parameters.

Before passing to the derivation of these results, we describe our
notation.
For an abelian variable $z\in\C$ and a base $p=e^{2\pi i\tau},$
Im$(\tau)>0$ (i.e., $|p|<1$), we define the infinite product
$$
(z;p)_\infty:=\prod_{j=0}^\infty(1-zp^j).
$$
Elliptic theta functions with characteristics are defined by the series
$$
\theta_{ab}(u)=\sum_{k\in \Z}e^{\pi i\tau(k+a/2)^2}
e^{2\pi i(k+a/2)(u+b/2)},
$$
where the variables $a$ and $b$ take values $0$ or $1$ and $u\in\C$.
The standard Jacobi theta functions are defined as
\begin{eqnarray*}
&& \theta_1(u|\tau)=\theta_1(u)=-\theta_{11}(u),
\\ &&
\theta_2(u|\tau)=\theta_2(u)=\theta_{10}(u)=\theta_1(u+1/2),
\\ &&
\theta_3(u|\tau)=\theta_3(u)=\theta_{00}(u)
=e^{\pi i \tau/4+\pi i u}\theta_1(u+1/2+\tau/2),
\\ &&
\theta_4(u|\tau)=\theta_4(u)=\theta_{01}(u)
=-ie^{\pi i \tau/4+\pi i u}\theta_1(u+\tau/2).
\end{eqnarray*}
It is convenient to set $\theta_a(u_1,\ldots,u_k):=\theta_a(u_1)\cdots
\theta_a(u_k)$ and $\theta_a(x\pm y):=\theta_a(x+y,x-y)$.

The $\theta_1(u)$-function is odd, $\theta_1(-u)=-\theta_1(u),$
and obeys the quasiperiodicity properties:
$$
\theta_1(u+1)=-\theta_1(u),\qquad
\theta_1(u+\tau)=-e^{-\pi i\tau -2\pi iu}\theta_1(u).
$$
The short theta function
$$
\theta(z;p):=(z;p)_\infty (pz^{-1};p)_\infty, \quad z\in\C^*,
$$
is related to $\theta_1(u)$ by the Jacobi triple product identity
$$
\theta_1(u)=ip^{1/8}e^{-\pi iu}(p;p)_\infty\theta(e^{2\pi iu};p).
$$
We shall need the duplication formula for $\theta_1$-function
$$
\theta_1(2u)=\frac{ip^{1/8}}{(p;p)_\infty^3}
\theta_1\left(u,u+\frac{1}{2},u+\frac{\tau}{2},u-\frac{1+\tau}{2}\right)
$$
and the addition formula
\begin{equation}
\theta_1(u\pm x,v\pm y)-\theta_1(u\pm y,v\pm x)=\theta_1(x\pm y,u\pm v)
\label{add-add}\end{equation}
or
\begin{equation}
\theta(ux^\pm,vy^\pm;p)-\theta(uy^\pm,vx^\pm;p)=\frac{v}{x}\theta(xy^\pm,uv^\pm;p),
\label{add-mult}\end{equation}
where $\theta(ux^\pm;p):=\theta(ux,ux^{-1};p):=\theta(ux;p)\theta(ux^{-1};p)$.

The standard elliptic gamma function has the form
\begin{equation}
 \Gamma_{\! p,q}(z) =\prod_{j,k=0}^\infty \frac{1-z^{-1}p^{j+1}q^{k+1}}
{1-zp^jq^k}, \quad |p|, |q|<1.
\label{ell-gamma}\end{equation}
We set
\begin{eqnarray*}
&& \Gamma_{\! p,q}(t_1,\ldots,t_k):=\Gamma_{\! p,q}(t_1)\cdots\Gamma_{\! p,q}(t_k),\quad
\\ &&
\Gamma_{\! p,q}(tz^\pm):=\Gamma_{\! p,q}(tz)\Gamma_{\! p,q}(tz^{-1}),\quad
\Gamma_{\! p,q}(z^{\pm2}):=\Gamma_{\! p,q}(z^2)\Gamma_{\! p,q}(z^{-2}).
\end{eqnarray*}
This function satisfies two finite difference equations of the first order
$$
\Gamma_{\! p,q}(qz)=\theta(z;p)\Gamma_{\! p,q}(z),\quad
\Gamma_{\! p,q}(pz)=\theta(z;q)\Gamma_{\! p,q}(z),
$$
the reflection equation
$$
\Gamma_{\! p,q}(z)\Gamma_{\! p,q}(pq/z)=1,
$$
the duplication formula
$$
\Gamma_{\! p,q}(z^2)=\Gamma_{\! p,q}(z,-z,q^{1/2}z,-q^{1/2}z, p^{1/2}z,-p^{1/2}z,
(pq)^{1/2}z,-(pq)^{1/2}z),
$$
and the limiting relation $\lim_{p\to0}\Gamma_{\! p,q}(z)=1/(z;q)_\infty$.

With the help of three pairwise incommensurate parameters
$\omega_{1,2,3}\in\C$, we define three base variables
\begin{eqnarray*}
&& q= e^{2\pi i\frac{\omega_1}{\omega_2}}, \quad
p=e^{2\pi i\frac{\omega_3}{\omega_2}}, \quad  r=e^{2\pi i\frac{\omega_3}{\omega_1}},
\end{eqnarray*}
and their particular modular transforms
\begin{eqnarray*}
\tilde q= e^{-2\pi i\frac{\omega_2}{\omega_1}}, \quad
\tilde p=e^{-2\pi i\frac{\omega_2}{\omega_3}},   \quad
\tilde r=e^{-2\pi i\frac{\omega_1}{\omega_3}}.
\end{eqnarray*}
The modified elliptic gamma function, which remains well defined for
$\omega_1/\omega_2>0$ (i.e., $|q|=1$),  has the form
\be
G(u;\mathbf{\omega})=
\Gamma_{\! p,q}(e^{2\pi i \frac{u}{\omega_2}})
\Gamma_{\tilde q,r}(re^{-2\pi i \frac{u}{\omega_1}})
= e^{-\pi i P(u)}\Gamma(e^{-2\pi i \frac{u}{\omega_3}};\tilde r,\tilde p),
\lab{unit-e-gamma}\ee
where $|p|, |r|<1$ and
$$
P\left(u+\sum_{m=1}^3\frac{\omega_m}{2}\right)
=\frac{u(u^2-\frac{1}{4}\sum_{m=1}^3\omega_m^2)}{3\omega_1\omega_2\omega_3}.
$$
It is  fixed uniquely by the following three equations
\begin{eqnarray}
&& f(u+\omega_1)=\theta(e^{2\pi i u/\omega_2};p)f(u),
\quad f(u+\omega_2) =\theta(e^{2\pi i u/\omega_1};r) f(u),
\nonumber \\ &&
 f(u+\omega_3) =e^{-\pi iB_{2,2}(u;\omega_1,\omega_2)} f(u),
\lab{e-gamma-eq}\\ &&
B_{2,2}(u;\omega_1,\omega_2)=\frac{u^2}{\omega_1\omega_2}
-\frac{u}{\omega_1}-\frac{u}{\omega_2}+
\frac{\omega_1}{6\omega_2}+\frac{\omega_2}{6\omega_1}+\frac{1}{2}
\nonumber \end{eqnarray}
and the normalization $f(\sum_{m=1}^3\omega_m/2)=1$.
The reflection equation has the form $G(a;\mathbf{\omega})G(b;\mathbf{\omega})=1$,
$a+b=\sum_{k=1}^3\omega_k$.
A more detailed description of the elliptic gamma functions can be
found in reviews \cite{rui:rev} or \cite{spi:thesis} and in the references therein.

\section{Continuous biorthogonal functions}

We investigate  elliptic hypergeometric integrals of the form
\begin{eqnarray}
&& I^{(m)}(t_1,\ldots,t_{2m+6})
=\kappa\int_\T\frac{\prod_{j=1}^{2m+6}\Gamma_{\! p,q}(t_jz^\pm)}
{\Gamma_{\! p,q}(z^{\pm2})}
\frac{dz}{z}, \quad
\label{ell-hyp-int}
\\ && \makebox[2em]{}
 \prod_{j=1}^{2m+6}t_j=(pq)^{m+1}, \quad
\kappa=\frac{(p;p)_\infty(q;q)_\infty}{4\pi i},
\nonumber\end{eqnarray}
where $\T$ is the unit circle and $|t_j|<1$.
For $m=0$, we obtain the elliptic beta integral \cite{spi:umn}
$$
I^{(0)}(t_1,\ldots,t_{6})=\prod_{1\leq j<k\leq 6}\Gamma_{\! p,q}(t_jt_k).
$$
Using this formula, it is easy to verify the following recursion relation
by mere change of the order of integrations
\ba \nonumber
&& I^{(m+1)}(t_1,\ldots,t_{2m+8})=\frac{\prod_{2m+5\leq k<l\leq 2m+8}
\Gamma_{\! p,q}(t_kt_l)}{\Gamma_{\! p,q}(\ve_m^2)}
\\ && \makebox[1em]{} \times
\kappa\int_\T\frac{\prod_{k=2m+5}^{2m+8}\Gamma_{\! p,q}(\ve_m^{-1}t_kw^\pm)}
{\Gamma_{\! p,q}(w^{\pm2})}
I^{(m)}(t_1,\ldots,t_{2m+4},\ve_m w,\ve_m w^{-1})\frac{dw}{w},
\label{key-rec}\ea
where $\ve_m=\sqrt{\prod_{k=2m+5}^{2m+8}t_k/pq}$.
It leads to the $m$-tuple integral representation for
$I^{(m)}$ similar to the Euler type integral representation for
the general plain hypergeometric function $_{m+1}F_m$,
i.e. $I^{(m)}$ can be interpreted as an elliptic analogue of the
$_{m+1}F_m$-function.
Recurrence \eqref{key-rec} is a special realization of an integral
analogue of the Bailey chains discovered in \cite{spi:bai}
(the Bailey chains technique is well known as a simplest tool for proving
the Rogers-Ramanujan type identities \cite{aar}).

For $m=0$, the integrand on the right-hand side of relation \eqref{key-rec}
contains the elliptic beta integral, computation of which yields \cite{spi:aa}
\be
V(t_1,\ldots,t_{8})=\prod_{1\leq k<l\leq 4}
\Gamma_{\! p,q}(t_kt_l,t_{k+4}t_{l+4})
V(\ve t_1,\ldots,\ve t_4, \ve^{-1}t_5,\ldots,\ve^{-1}t_8),
\label{E7}\ee
where $V(\underline{t})=I^{(1)}(\underline{t})$ is an elliptic
analogue of the Gauss hypergeometric function
and $\ve=\sqrt{t_5t_6t_7t_8/pq}=\sqrt{pq/t_1t_2t_3t_4}$.
This is the key reflection extending the obvious $S_8$-group
of symmetries in parameters of $V(\underline{t})$
to the Weyl group for the exceptional root system $E_7$.

For $m=1$, we obtain the equality
\ba \nonumber
&& I^{(2)}(t_1,\ldots,t_{10})=\frac{\prod_{7\leq k<l\leq 10}
\Gamma_{\! p,q}(t_kt_l)}{\Gamma_{\! p,q}(\ve_1^2)}
\\ && \makebox[4em]{} \times
\kappa\int_\T\frac{\prod_{k=7}^{10}\Gamma_{\! p,q}(\ve_1^{-1}t_kw^\pm)}
{\Gamma_{\! p,q}(w^{\pm2})}
V(t_1,\ldots,t_6,\ve_1 w,\ve_1 w^{-1})\frac{dw}{w},
\label{m=1}\ea
where $\ve_1=\sqrt{t_7t_8t_9t_{10}/pq}$.
After fixing parameters $t_5t_7=t_6t_{8}=pq$, the integral in
the left hand side gets reduced to the elliptic beta integral
and, so, it can be computed explicitly. This yields the relation
\ba\nonumber
&& \prod_{1\leq k<l\leq 4}\Gamma_{\! p,q}(t_kt_l)\prod_{k=1}^4\Gamma_{\! p,q}(t_kt_9,t_kt_{10})
 =
\frac{\prod_{7\leq k<l\leq 10}\Gamma_{\! p,q}(t_kt_l)}
{\Gamma_{\! p,q}(t_9t_{10})\Gamma_{\! p,q}(\ve_1^2)}
\\ && \makebox[2em]{} \times
\kappa\int_\T \frac{\prod_{j=7}^{10}\Gamma_{\! p,q}(\ve_1^{-1}t_j w^\pm)}
{\Gamma_{\! p,q}(w^{\pm2})}
V(t_1,\ldots,t_4,pq/t_7,pq/t_8, \ve_1w,\ve_1w^{-1})\frac{dw}{w},
\nonumber\ea
where $t_1t_2t_3t_4t_9t_{10}=pq$.
We denote now $t_9=s\xi, t_{10}=s\xi^{-1}$ and obtain
\ba\nonumber
&& \prod_{1\leq k<l\leq 4}\Gamma_{\! p,q}(t_kt_l)\prod_{j=1}^4\Gamma_{\! p,q}(t_js\xi^\pm)
=\frac{\Gamma_{\! p,q}(t_7s\xi^\pm,t_{8}s\xi^\pm,t_7t_{8})}{\Gamma_{\! p,q}(\ve_1^2)}
\\ && \makebox[4em]{} \times
\kappa\int_\T \frac{\Gamma_{\! p,q}(
\ve_1^{-1}t_7w^\pm,\ve_1^{-1}t_{8}w^\pm,\ve_1^{-1}s\xi^\pm w^\pm)}
{\Gamma_{\! p,q}(w^{\pm2})}
 \nonumber \\&& \makebox[6em]{} \times
V(t_1,\ldots,t_4, pq/t_7,pq/t_8,\ve_1w,\ve_1w^{-1})\frac{dw}{w},
\label{key-rel'}\ea
where $t_1t_2t_3t_4s^2=pq$ and $\ve_1=s\sqrt{t_7t_{8}/pq}$. After the replacements
$$
 t_1=c,\quad t_2=d, \quad
t_3=s^{-1}\sqrt{\frac{pq}{cd}}\,x,\quad t_4=s^{-1}\sqrt{\frac{pq}{cd}}\, x^{-1},
\quad t_7=a,\quad t_{8}=b,
$$
the condition $t_1t_2t_3t_4s^2=pq$ is satisfied automatically
and $\ve_1=s\sqrt{ab/pq}$. Introducing the basis vectors
\be
\phi(w;a,b|\xi;s)=\Gamma_{\! p,q}(sa\xi^\pm,sb\xi^\pm, \sqrt{\frac{pq}{ab}}w^\pm\xi^\pm),
\lab{phi}\end{equation}
we can  rewrite the key relation  \eqref{key-rel'} as
\be
\phi(x;c,d|\xi;s)=\kappa\int_\T R(c,d,a,b;x,w|s)\phi(w;a,b|\xi;s)\frac{dw}{w},
\lab{key-rel-fin}\ee
where, after the substitution $s=\sqrt{pq/\rho}$,
\begin{eqnarray}\nonumber
&& R(c,d,a,b;x,w|s)=\frac{\Gamma_{\! p,q}(ab,
\sqrt{\frac{a\rho}{b}}w^\pm,\sqrt{\frac{b\rho}{a}}w^\pm)}
{\Gamma_{\! p,q}(cd,\sqrt{\frac{c\rho}{d}}x^\pm,
\sqrt{\frac{d\rho}{c}}x^\pm)}\frac{1}{\Gamma_{\! p,q}(\frac{ab}{\rho},
\frac{\rho}{cd},w^{\pm2})}
\\ && \makebox[4em]{} \times
V\left(c,d,\sqrt{\frac{\rho}{cd}}x,
\sqrt{\frac{\rho}{cd}}x^{-1},\frac{pq}{a},\frac{pq}{b},
\sqrt{\frac{ab}{\rho}}w,\sqrt{\frac{ab}{\rho}}w^{-1}\right).
\label{R}\end{eqnarray}
Using the transformation \eqref{E7}, we can bring this $R$-function to the form
\begin{eqnarray*}\nonumber
&& R(c,d,a,b;x,w|s)=
\frac{1}{\Gamma_{\! p,q}(\frac{pq}{ab},\frac{ab}{pq},w^{\pm2})}
\\ && \makebox[4em]{} \times
V\left(s c,s d,\sqrt{\frac{pq}{cd}}x,
\ve\sqrt{\frac{pq}{cd}}x^{-1},\frac{pq}{as},\frac{pq}{bs},
\sqrt{\frac{ab}{pq}}w,\sqrt{\frac{ab}{pq}}w^{-1}\right).
\end{eqnarray*}
The kernel $R$ does not depend on the variable $\xi$. Relation
\eqref{key-rel-fin} looks therefore like a rotation in some
functional space spanned by the basis vectors $\phi$, which is
performed by the ``matrix" $R$ with the continuous indices $x$ and $w$.

We can denote $\xi=e^{i\chi},\,w=e^{i\theta},$ and $x=e^{i\varphi}$
and let $\chi, \theta, $ and $\varphi$ take real values only.
Then, because of the symmetries
$$
\phi(w;a,b|\xi;s)=\phi(w^{-1};a,b|\xi;s)=\phi(w;a,b|\xi^{-1};s),
$$
we see that this function depends actually on the variables
$y:=\cos\theta$ and $\cos\chi$. Similarly, the function
$$
R(a,b,c,d;x,w|s)=R(a,b,c,d;x^{-1},w|s)=R(a,b,c,d;x,w^{-1}|s)
$$
depends on $y$ and $\cos\varphi$. Using the elementary relation
$$
\int_\T f(\cos\theta)\frac{dw}{iw}=\int_0^{2\pi}f(\cos\theta)d\theta
=2\int_{-1}^1f(y)\frac{dy}{\sqrt{1-y^2}},
$$
we can write
\be
\phi(x;c,d|\xi;s)=2i\kappa
\int_{-1}^1 R(c,d,a,b;x,e^{i\theta}|s)
\phi(e^{i\theta};a,b|\xi;s) \frac{dy}{\sqrt{1-y^2}}.
\lab{key-real}\ee
Taking the limits $c\to a$ and $d\to b$ and denoting $v=\cos\varphi$,
we obtain the relation
$$
\lim_{c\to a,d\to b}
R(c,d,a,b;e^{i\varphi},e^{i\theta})
=\frac{2\pi\sqrt{1-y^2}}{(p;p)_\infty(q;q)_\infty}\,\delta(v-y)
$$
in the distributional sense. Formally, this equality can be
rewritten as
\begin{eqnarray}\nonumber
&& V\left(a,b,\sqrt{\frac{\rho}{ab}}e^{i\varphi},
\sqrt{\frac{\rho}{ab}}e^{-i\varphi},\frac{pq}{a},\frac{pq}{b},
\sqrt{\frac{ab}{\rho}}e^{i\theta},\sqrt{\frac{ab}{\rho}}e^{-i\theta}\right)
\\ && \makebox[4em]{}
=\Gamma_{\! p,q}\left(\frac{ab}{\rho},\frac{\rho}{ab},e^{\pm 2i\theta}\right)
\frac{2\pi\sqrt{1-y^2}}{(p;p)_\infty(q;q)_\infty}\,\delta(v-y).
\label{V-delta}\end{eqnarray}

Using twice equality \eqref{key-rel-fin} with different choices of
parameters on the right-hand side, we obtain immediately the
self-reproducing property for the kernel
\begin{equation}
\kappa\int_\T R(a,b,c,d;x,w|s)R(c,d,e,f;w,z|s)\frac{dw}{w}
=R(a,b,e,f;x,z|s)
\label{repr}\end{equation}
and the biorthogonality relation
\begin{eqnarray}\label{biort}
&&
\int_{-1}^1 R(a,b,c,d;e^{i\varphi},e^{i\theta}|s)
R(c,d,a,b;e^{i\theta},e^{i\varphi'}|s)
\frac{dy}{\sqrt{1-y^2}}=\frac{\sqrt{1-v^2}}{(2i\kappa)^2}\,\delta(v-v'),
\nonumber\end{eqnarray}
where $v=\cos\varphi$ and $v'=\cos\varphi'$ can be considered as continuous spectral
variables. Substituting the explicit expressions for the $R$-functions,
we rewrite the biorthogonality relation as
\begin{eqnarray}\nonumber
&&
\int_{-1}^1 \frac{1}{\Gamma_{\! p,q}(e^{\pm2i\theta})}
V\left(a,b,\sqrt{\frac{\rho}{ab}}e^{i\varphi},
\sqrt{\frac{\rho}{ab}}e^{-i\varphi},\frac{pq}{c},\frac{pq}{d},
\sqrt{\frac{cd}{\rho}}e^{i\theta},\sqrt{\frac{cd}{\rho}}e^{-i\theta}\right)
\\ && \makebox[2em]{} \times
V\left(c,d,\sqrt{\frac{\rho}{cd}}e^{i\theta},
\sqrt{\frac{\rho}{cd}}e^{-i\theta},\frac{pq}{a},\frac{pq}{b},
\sqrt{\frac{ab}{\rho}}e^{i\varphi'},\sqrt{\frac{ab}{\rho}}e^{-i\varphi'}\right)
\frac{dy}{\sqrt{1-y^2}}
\nonumber \\ && \makebox[2em]{}
=\Gamma_{\! p,q}\left(\frac{ab}{\rho},\frac{\rho}{ab},
\frac{cd}{\rho},\frac{\rho}{cd},e^{\pm2i\varphi}\right)
\frac{\sqrt{1-v^2}}{(2i\kappa)^2}\,\delta(v-v').
\label{V-ort}\end{eqnarray}

\section{A scalar product and biorthogonality of the basis vectors}

It is necessary to establish now the Hilbert space content of the construction
described above. It is desirable to connect the $V$-function with a
scalar product of some basis vectors of this space.

Multiplying equality \eqref{key-rel-fin} by $\phi(z;pq/e,pq/f|\xi;s^{-1})$
and integrating it over $\xi$
with the measure $\kappa\int_\T \Gamma_{\! p,q}^{-1}(\xi^{\pm2})
d\xi/\xi$, we obtain on the left-hand side
\begin{equation}
\kappa\int_\T \frac{\phi(z;\frac{pq}{e},\frac{pq}{f}|\xi;s^{-1})\phi(x;c,d|\xi;s)}
{\Gamma_{\! p,q}(\xi^{\pm2})}\frac{d\xi}{\xi}
=V\left(st_1,\ldots,st_4,s^{-1}t_5,\ldots,s^{-1}t_8\right),
\label{prel}\end{equation}
where $t_1=c,t_2=d,t_{3,4}=x^{\pm1}s^{-1}\sqrt{pq/cd}$ and
$t_5=pq/e,t_6=pq/f,t_{7,8}= z^{\pm1}s\sqrt{ef/pq}$. On the
right-hand side, we find
$$
\kappa\int_\T  R(c,d,a,b;x,w|s)
V\left(st_1',\ldots,st_4',s^{-1}t_5,\ldots,s^{-1}t_8\right)\frac{dw}{w},
$$
where $t_1'=a,t_2'=b,t_{3,4}'=w^{\pm1}s^{-1}\sqrt{pq/ab}$.
If we apply $E_7$-transformation \eqref{E7}
to $V$-functions on both sides of the equality, we obtain relation \eqref{repr},
which shows the self-consistency of the consideration.

Taking the limits $e\to c$, $f\to d$, and applying equality \eqref{V-delta}
to the right-hand side of \eqref{prel}, we find the biorthogonality relation
\begin{eqnarray}\nonumber
&& \kappa\int_\T \frac{\phi(e^{i\varphi'};\frac{pq}{c},\frac{pq}{d}|\xi;s^{-1})
\phi(e^{i\varphi};c,d|\xi;s)}
{\Gamma_{\! p,q}(\xi^{\pm2})}\frac{d\xi}{\xi}
\\ && \makebox[4em]{}
=\Gamma_{\! p,q}\left(\frac{pq}{cd},\frac{cd}{pq},e^{\pm2i\varphi}
\right)\sqrt{1-v^2}\, \delta(v-v').
\label{phi-bio}\end{eqnarray}

The latter result inspires introduction of the following
inner product. For two functions of $z\in\C$ with the
property $\psi(z)=\psi(z^{-1})$, we define
\begin{equation}
\langle \chi(z),\psi(z)\rangle =\kappa\int_{\T}\frac{\chi(z)\psi(z)}
{\Gamma_{\! p,q}(z^{\pm 2})}
\frac{dz}{z}.
\label{sp}\end{equation}
From the reflection equation for the elliptic gamma function, we find
for $z\in\T$ and real $(p, q)$-variables
$$
\frac{1}{\Gamma_{\! p,q}(z^{\pm 2})}=\theta(z^2;p)\theta(z^{-2};q)
=\left|(1-z^2)(z^2;p)_\infty(qz^2;q)_\infty\right|^2,
$$
i.e. our measure is positive definite. Similarly, for $p=q^*$ we obtain
a  positive measure as well
$$
\frac{1}{\Gamma_{\! p,q}(z^{\pm 2})}=|\theta(z^2;p)|^2
=\left|(z^2;p)_\infty(qz^2;q)_\infty\right|^2.
$$

We introduce now the finite difference operator \cite{rai:trans,ros:elementary}
$$
D_R(a,b,c,d;p;q):=\frac{\theta(az,bz,cz,dz;p)}{\theta(z^2;p)}T_{z,q}^{1/2}
+\frac{\theta(az^{-1},bz^{-1},cz^{-1},dz^{-1};p)}{\theta(z^{-2};p)}T_{z,q}^{-1/2},
$$
where  $T_{z,q}^{\pm1/2}f(z)=f(q^{\pm1/2}z)$ is a $q$-shift operator.
The operator conjugated  to $D_R$ is defined by the equality
$$
\langle \chi,D_R(a,b,c,d;p;q)\psi\rangle=\langle D_R^*(a,b,c,d;p;q)\chi,\psi\rangle.
$$
Assuming that functions $\psi(z)$ and $\chi(z)$ do not have
singularities in the annulus $|q|^{1/2}<|z|<|q|^{-1/2}$, we find
$$
D_R^*(a,b,c,d;p;q)=D_R\left(\frac{pq^{1/2}}{a},\frac{pq^{1/2}}{b},
\frac{pq^{1/2}}{c}, \frac{pq^{1/2}}{d};p;q\right).
$$
Clearly, $D_R^{**}=D_R$.
For further considerations, it is convenient to scale one of the parameters
by $p$ and define $D(a,b,c,d;p;q):=-aD_R(pa,b,c,d;p;q)$, so that
\begin{equation}
D(a,b,c,d;p;q)=
\frac{\theta(az,bz,cz,dz;p)}{z\theta(z^2;p)}T_{z,q}^{1/2}
+\frac{\theta(az^{-1},bz^{-1},cz^{-1},dz^{-1};p)}
{z^{-1}\theta(z^{-2};p)}T_{z,q}^{-1/2}
\label{D-oper}\end{equation}
and
$$
D^*(a,b,c,d;p;q)=\frac{cd}{q^{1/2}}D\left(\frac{pq^{1/2}}{a},\frac{pq^{1/2}}{b},
\frac{q^{1/2}}{c}, \frac{q^{1/2}}{d};p;q\right).
$$

We look now for solutions of the generalized eigenvalue problem
for $D$-operators with the fixed product $\rho:=abcd$,
\begin{equation}
D(a,b,c,d;p;q)f(z)=\lambda D(a,b,c',d';p;q)f(z),
\label{GEVP}\end{equation}
where $cd=c'd'$.
We parametrize the spectral variable $\lambda$
by an elliptic function of the second order
$$
\lambda\equiv \lambda(w)=\frac{\theta(w\sqrt{c/d},w\sqrt{d/c};p)}
{\theta(w\sqrt{c'/d'},w\sqrt{d'/c'};p)},
$$
where $w\in\C$ is an arbitrary parameter, and indicate in $f(z)$
dependence on all parameters
$$
f(z):=f(z;w;q^{1/2}a,q^{1/2}b;\rho),
$$
where the parameters $a$ and $b$ are scaled by $q^{1/2}$ for a later convenience.
 Equation \eqref{GEVP} can be rewritten then as
\begin{eqnarray*}
&& \makebox[-2em]{}
 \frac{f(q^{1/2}z;w;q^{1/2}a,q^{1/2}b;\rho)}
{f(q^{-1/2}z;w;q^{1/2}a,q^{1/2}b;\rho)}
=-z^2\frac{\theta(z^2,az^{-1},bz^{-1};p)}
{\theta(z^{-2},az,bz;p)}
\\ &&  \makebox[-1em]{}
\times
\frac{\theta(cz^{-1},dz^{-1},w\sqrt{c'/d'},w\sqrt{d'/c'};p)
- \theta(c'z^{-1},d'z^{-1},w\sqrt{c/d},w\sqrt{d/c};p)}
{\theta(cz,dz,w\sqrt{c'/d'},w\sqrt{d'/c'};p)
- \theta(c'z,d'z,w\sqrt{c/d},w\sqrt{d/c};p)},
\end{eqnarray*}
which, after application of the addition formula for $\theta(z;p)$-functions,
reduces to
\be
\frac{f(q^{1/2}z;w;q^{1/2}a,q^{1/2}b;\rho)}
{f(q^{-1/2}z;w;q^{1/2}a,q^{1/2}b;\rho)}=z^4
\frac{\theta(az^{-1},bz^{-1},\sqrt{cd}w^\pm z^{-1};p)}
{\theta(az,bz,\sqrt{cd} w^\pm z;p)}.
\lab{f-eq}\ee
A particular solution of this equation valid for $|q|<1$
has the form
\begin{equation}
f(z;w;a,b;\rho)=\Gamma_{\! p,q}\left(\frac{pq}{a} z^\pm,
\frac{pq}{b} z^\pm, \sqrt{\frac{ab}{\rho}}w^\pm z^\pm\right).
\label{basis}\end{equation}
The general solution at $|q|<1$ is obtained from  \eqref{basis} after
multiplying it by arbitrary function $\varphi(z)$ satisfying the equality
$\varphi(qz)=\varphi(z)$. The case $|q|=1$ will be considered below,
but the regime $|q|>1$ will be omitted for brevity, since it
turns out to be equivalent to the case $|q|<1$ after a change
of parameters.

Our particular solution depends only on three parameters
$a,b,$ and $\rho$, i.e. two variables, say $c$ and $c'$,
drop out being plain gauge parameters. A special choice of parameters
reduces $f(z;w;a,b;\rho)$ to products of theta functions
\begin{equation}
f\left(z;q^k\sqrt{\frac{a}{bq^{N}}};a,b;q^{-N}\right)
=\prod_{j=0}^{k-1}\theta(q^jaz^\pm;p)\prod_{j=0}^{N-k-1}
\theta(q^jbz^\pm;p),\quad  k,N\in\N,
\label{dis-bas}\end{equation}
which coincide with the intertwining vectors of Takebe \cite{tak}
and the basis vectors used by Rosengren in \cite{ros:elementary,ros:sklyanin}.

Since the function $f(z;w;a,b;\rho)$ is symmetric in $p$ and $q$,
we can demand  it satisfy one more generalized
eigenvalue problem obtained from \eqref{gevp2} by permutation of $p$ and $q$:
\begin{eqnarray}\nonumber
&& D(a,b,c,d;q;p)f(z;w;p^{1/2}a,p^{1/2}b;\rho)
\\ && \makebox[2em]{}
=\frac{\theta(w\sqrt{c/d},w\sqrt{d/c};q)}
{\theta(w\sqrt{c'/d'},w\sqrt{d'/c'};q)}
D(a,b,c',d';q;p)f(z;w;p^{1/2}a,p^{1/2}b;\rho),
\label{GEVP-perm}\end{eqnarray}
where $cd=c'd'$.  For incommensurate $p$ and $q$ (which we assume
throughout the paper), this requirement strongly restricts the space of
solutions of our eigenvalue problem
since the equations $\varphi(qz)=\varphi(pz)=\varphi(z)$
are satisfied only by a constant function. Our basis vectors
$f(z;w;a,b;\rho)$ become thus fixed uniquely up to a constant multiplier
which we set equal to 1 for concreteness.

Function \eqref{basis} satisfies another type of the generalized
eigenvalue problem
\be
D(a,b,c,d;p;q)f(z;w;q^{1/2}a,q^{1/2}b;\rho)=d^{-1}\theta(w^\pm\sqrt{d/c},
cd;p)f(z;w;a,b;\rho).
\label{gevp2}\ee
Indeed, dividing this relation by $f(z;w;a,b)$,
we obtain on the left-hand side
\begin{eqnarray*}
&& \frac{\theta(az,bz,cz,dz;p)}{z\theta(z^2;p)}
\frac{f(q^{1/2}z;w;q^{1/2}a,q^{1/2}b;\rho)}{f(z;w;a,b;\rho)}
\\ && \makebox[2em]{}
+\frac{\theta(az^{-1},bz^{-1},cz^{-1},dz^{-1};p)}{z^{-1}\theta(z^{-2};p)}
\frac{f(q^{-1/2}z;w;q^{1/2}a,q^{1/2}b;\rho)}{f(z;w;a,b;\rho)}.
\end{eqnarray*}
After substitution of the explicit form of $f(z;w;a,b;\rho)$, this expression
reduces to
$$
\frac{\theta(cz,dz,w^\pm z\sqrt{ab/\rho};p)
-\theta(c^{-1}z,d^{-1}z,w^\pm z\sqrt{\rho/ab};p)}
{z\theta(z^2;p)}.
$$
Addition formula \eqref{add-mult} with $u=z\sqrt{cd},\, x=\sqrt{c/d},\,
v=z\sqrt{ab/\rho},\, y=w$ yields an expression required for the
right-hand side in \eqref{gevp2}. Equality \eqref{GEVP}
is easily deduced from \eqref{gevp2} by noticing that $f(z;w;a,b;\rho)$
does not change after replacements of $c$ and $d$ by $c'$ and $d'$
with the constraint $cd=c'd'$.

The following commutation relation \cite{rai:abelian}
\begin{eqnarray}\nonumber
&& D(a,b,c',d';p;q)D(q^{1/2}a,q^{1/2}b,q^{-1/2}c,q^{-1/2}d;p;q)
\\ && \quad
=D(a,b,c,d;p;q)D(q^{1/2}a,q^{1/2}b,q^{-1/2}c',q^{-1/2}d';p;q),
\lab{ccr}\end{eqnarray}
where $cd=c'd'$, is a simple consequence of \eqref{gevp2}.
Indeed, the action of both sides of this equality on the function
$f(z;w;qa,qb;\rho)$ yields the same result. The straightforward verification
of the identity is also possible (coincidence of the terms in front of
the operators $T^{\pm1}_{z,q}$ is evident, and the scalar terms coincide
after application of the addition formula for theta functions).

Consider the generalized eigenvalue problem for the conjugated
operators
$$
D^*(a,b,c,d;p;q)g(z)=\lambda(v) D^*(a,b,c',d';p;q)g(z),
$$
where $cd=c'd'$ and $\lambda(v)$ is a new spectral variable.
Evidently, this equation is obtained from \eqref{GEVP}
after replacements $a\to \tilde a= pq^{1/2}/a, b\to \tilde b= pq^{1/2}/b,
c\to \tilde c= q^{1/2}/c, $ $d\to \tilde d= q^{1/2}/d$, and
$w\to v$. We find thus its particular solution
\begin{equation}
g(z):= g(z;v;a,b;\rho)= \Gamma_{\! p,q}\left(a z^\pm,
b z^\pm, \sqrt{\frac{\rho}{ab}}v^\pm z^\pm\right).
\label{basis-dual}\end{equation}
The biorthogonality relation of interest emerges as follows
\begin{eqnarray}\nonumber
&&
0=\langle \left(D^*(a,b,c,d;p;q)-\lambda(v)D^*(a,b,c',d';p;q)\right)g(z;v;a,b;\rho),
\\ && \makebox[2em]{}
 f(z;w;q^{1/2}a,q^{1/2}b;\rho)\rangle
=\langle g(z;v;a,b;\rho),
\nonumber \\ && \makebox[2em]{}
\left(D(a,b,c,d;p;q)-\lambda(v)D(a,b,c',d';p;q)\right)
 f(z;w;q^{1/2}a,q^{1/2};\rho)\rangle
\nonumber\\ && \makebox[2em]{}
=\left(1-\frac{\lambda(v)}{\lambda(w)}\right)\langle g(z;v;a,b;\rho),
D(a,b,c,d;p;q)f(z;w;q^{1/2}a,q^{1/2}b;\rho)\rangle
\nonumber\\ && \makebox[2em]{}
=\left(1-\frac{\lambda(v)}{\lambda(w)}\right)d^{-1}\theta(w^\pm\sqrt{d/c},cd;p)
\langle g(z;v;a,b;\rho),f(z;w;a,b;\rho)\rangle.
\label{GEVP-dual}\end{eqnarray}
For $\lambda(v)\neq \lambda(w)$, we find thus
$\langle g(z;v;a,b;\rho),f(z;w;a,b;\rho)\rangle=0,$ which
essentially coincides with relation \eqref{V-delta}.
Restrictions on the parameters needed for conjugation of the operator
$D$ without crossing singularities are $|p|<|a|, |b|<|q|^{1/2}$
and $|\rho/abq|<|v|^2< |abq/\rho|$, $|ab/\rho q |<|w|^2<|\rho q/ab|$.

To conclude, solutions of the generalized eigenvalue problems for the $D$-operator
and its conjugate define a basis rotated by the $R$-function
and its dual. An overlap of these bases defines the $V$-function:
\begin{eqnarray} \makebox[-1em]{}
\kappa\int_\T\frac{g(z;v;a,b;e)f(z;w;c,d;e)}
{\Gamma_{\! p,q}(z^{\pm2})}\frac{dz}{z}
= V\left(a,b,\sqrt{\frac{e}{ab}}v^\pm,
\frac{pq}{c},\frac{pq}{d},
\sqrt{\frac{cd}{e}}w^\pm
\right),
\label{V-fact}\end{eqnarray}
where $V(\ldots,ax^\pm,\ldots)=V(\ldots,ax,ax^{-1},\ldots)$.

We can rewrite also the key relation  \eqref{key-rel-fin} in a compact form
\begin{eqnarray}\nonumber
&& \makebox[-1em]{}
\Gamma_{\! p,q}(\alpha z^\pm\xi^\pm,\beta x^\pm\xi^\pm)
=\kappa\int_\T r(\alpha,\beta,\gamma,\delta;z,x;t,w)
\Gamma_{\! p,q}(\gamma t^\pm\xi^\pm,\delta w^\pm\xi^\pm)\frac{dw}{w},
\\ &&  \makebox[-1em]{}
r(\alpha,\beta,\gamma,\delta;z,x;t,w)=\frac{1}
{\Gamma_{\! p,q}(\delta^{\pm2},w^{\pm2})}
V\left(\beta x^\pm,\alpha z^\pm,
\frac{pq}{\gamma} t^\pm,\frac{w^{\pm}}{\delta}\right),
\nonumber\end{eqnarray}
where $\alpha\beta=\gamma\delta$.

\section{Elliptic hypergeometric equation}

We recall the derivation of the elliptic hypergeometric equation
together with its general solution
\cite{spi:thesis,spi:cs}. Repeating transformation \eqref{E7}
with the parameters $\ve t_3,\ve t_4, \ve^{-1}t_5,
\ve^{-1}t_6$, playing the role of $t_1,\dots,t_4$, and permuting
$t_3,t_4$ with $t_5,t_6$ in the result, we find
\be
V(\underline{t})=\prod_{j,k=1}^4
\Gamma_{\! p,q}(t_jt_{k+4} )\ V(\sqrt{T}/t_1,\ldots,\sqrt{T}/t_4,
\sqrt{U}/t_5,\ldots,\sqrt{U}/t_8),
\lab{E7-2}\ee
where $ T=t_1t_2t_3t_4$ and $ U=t_5t_6t_7t_8$.
Equating right-hand sides of \eqref{E7} and \eqref{E7-2},
we obtain the third symmetry transformation
\be
V(\underline{s})=\prod_{1\le j<k\le 8}\Gamma_{\! p,q}(s_js_k )\,
V(\sqrt{pq}/s_1,\ldots,\sqrt{pq}/s_8).
\lab{E7-3}\ee
All these relations represent different Weyl group transformations
for the exceptional root system $E_7$ in the parameter space.

If we set $v=t_6, x=t_7, y=z,$ and $u=t_8$ in the
addition formula \eqref{add-mult}, then it becomes equivalent
to the equation
\begin{eqnarray}\nonumber &&
t_7\theta(t_8t_7^\pm;p)\Delta(z,t_1,\ldots,t_5,qt_6,t_7,t_8)
-t_7\theta(t_6t_7^\pm;p)\Delta(z,t_1,\ldots,t_7,qt_8)
 \\ && \makebox[4em]{}
=t_6\theta(t_8t_6^\pm;p)\Delta(z, t_1,\ldots,t_6, qt_7,t_8),
\nonumber \end{eqnarray}
where $\Delta(z,\underline{t})=\prod_{k=1}^8\Gamma_{\! p,q}(t_kz^\pm)/
\Gamma_{\! p,q}(z^{\pm2})$
is $V$-function's kernel. Integration over $z$ along an appropriate
contour $C$ yields
\begin{equation}
t_7\theta\left(t_8t_7^\pm;p\right)V(qt_6)
-t_6\theta\left(t_8t_6^\pm;p\right)V(qt_7)
=t_7\theta\left(t_6t_7^\pm;p\right) V(qt_8),
\label{c1}\end{equation}
where $V(qt_j)$ denote $V$-functions with the parameters
 $t_j$ replaced with $qt_j$ and the balancing condition
$\prod_{j=1}^8t_j=p^2q$.
Substituting transformation \eqref{E7-3} in this equation, we obtain
\begin{eqnarray}\nonumber
&& t_6\theta\left(\frac{t_7}{t_8};p\right)\prod_{k=1}^5
\theta\left(t_6t_k;p\right)V(qt_7,qt_8)
-t_7\theta\left(\frac{t_6}{t_8};p\right)\prod_{k=1}^5
\theta\left(t_7t_k;p\right)V(qt_6,qt_8)
\\ && \makebox[2em]{}
=t_6\theta\left(\frac{t_7}{t_6};p\right)
\prod_{k=1}^5\theta(t_8t_k;p) V(qt_6,qt_7),
\label{eq2}\end{eqnarray}
where $\prod_{j=1}^8t_j=p^2$.

We permute now parameters $t_5$ and $t_8$ in the latter equation,
change the parameters $t_5,t_7\to t_5/q,t_7/q$,
and exclude the function $V(qt_6,q^{-1}t_7)$ with the help of
equality \eqref{c1} where the parameter $t_7$ is replaced with $t_7/q$.
The resulting relation looks like
\begin{eqnarray}\nonumber
&& \frac{\theta(t_7t_6^\pm/q;p)}{\theta(t_8t_6^\pm;p)} V(q^{-1}t_7,qt_8)
-V(\underline{t})
\\ &&  \makebox[2em]{}
=\frac{\theta(t_7/qt_6,t_5t_8/q,qt_8/t_7;p)}
{\theta(t_8t_6^\pm,qt_6/t_5;p)}
\prod_{j=1}^4\frac{\theta(t_5t_j/q;p)}{\theta(t_7t_j/q;p)}V(q^{-1}t_5,qt_6)
\nonumber  \\ &&  \makebox[4em]{}
-\frac{\theta(qt_8/t_7,t_7/t_5;p)}
{\theta(t_8/t_6,qt_6/t_5;p)}\prod_{j=1}^4\frac{\theta(t_6t_j;p)}
{\theta(t_7t_j/q;p)} V(\underline{t}).
\label{key-cont}\end{eqnarray}
Multiplying it by the function
\begin{equation}
\mathcal{A}(t_1,\ldots,t_8;p)
=\frac{\theta(t_8/qt_6,t_6t_8,t_6/t_8;p)}{\theta(t_8/t_7,t_7/qt_8,t_7t_8/q;p)}
\prod_{j=1}^5\frac{\theta(t_7t_j/q;p)}{\theta(t_6t_j;p)}
\label{potential}\end{equation}
and symmetrizing in parameters $t_7$ and $t_8$
yield the elliptic hypergeometric equation
\begin{equation}
\mathcal{A}(t_1,\ldots,t_7,t_8;p)\left(\frac{\theta(t_7t_6^\pm/q;p)}{\theta(t_8t_6^\pm;p)} V(q^{-1}t_7,qt_8)
-V(\underline{t})\right) +(t_7\leftrightarrow t_8)+V(\underline{t})=0,
\label{eh}\end{equation}
which was derived in \cite{spi:thesis} in a slightly different way.

We set $t_7=cx^{-1},\, t_8=cx$ and denote
\begin{equation}
\mathcal{D}(t_1,\ldots,t_8)=a(t_1,\ldots,t_7,t_8)(T_{x,q}-1)+
a(t_1,\ldots,t_8,t_7)(T^{-1}_{x,q}-1)+\kappa(t_1,\ldots,t_8),
\label{D-eheq}\end{equation}
where $T_{x,q}f(x)=f(qx)$ and
\begin{eqnarray*}
&& a(t_1,\ldots,t_8)=\frac{\theta(t_8/qt_6,t_6t_8,t_6/t_8;p)}
                 {\theta(t_8/t_7,t_7/qt_8;p)}
\prod_{k=1}^5\theta(t_7t_k/q;p), \quad
\\ &&
\kappa(t_1,\ldots,t_8)=\theta(t_7t_8/q;p)\prod_{k=1}^5\theta(t_6t_k;p).
\end{eqnarray*}
Equation \eqref{eh} may be rewritten then as a linear second order
$q$-difference equation $\mathcal{D}(\underline{t})U(\underline{t})=0$,
where $U(\underline{t})=V(\underline{t})/\Gamma_{\! p,q}(t_7t_6^\pm,t_8t_6^\pm).$
Other linearly independent solutions of this equation can be
obtained by replacing $t_j\to pt_j$, for some fixed $j=1,\ldots,5$,
or $c\to pc$, or $x\to px$, etc \cite{spi:thesis}.

It is not difficult to verify the operator identity
\begin{equation}
\mathcal{D}(\underline{t}')-\frac{\theta(t_3/t_3',t_3/t_4';p)}
{\theta(t_3/t_3'',t_3/t_4'';p)}\mathcal{D}(\underline{t}'')
=\frac{t_3}{t_3'}\frac{\theta(t_3'/t_3'',t_3'/t_4'';p)}
{\theta(t_3/t_3'',t_3/t_4'';p)}\mathcal{D}(\underline{t}),
\label{op-ident}\end{equation}
where in the primed set of parameters $\underline{t}'$ variables $t_3$ and $t_4$
are replaced by arbitrary $t_3'$ and $t_4'$ with the constraint $t_3t_4=t_3't_4'$
 and similarly for the double primed parameters.
From the elliptic hypergeometric equation
$\mathcal{D}(\underline{t})U(\underline{t})=0$, we find that the ratio
$$
\frac{\mathcal{D}(\underline{t}') U(\underline{t})}{\theta(t_3/t_3',t_3/t_4';p)}=
\frac{\mathcal{D}(\underline{t}'') U(\underline{t})}{\theta(t_3/t_3'',t_3/t_4'';p)}
$$
is independent on the auxiliary parameter $t_3'$ (or $t_4'$).
To find its explicit form, we multiply contiguous relation \eqref{key-cont}
by $\mathcal{A}(\underline{t}')$ and symmetrize it in the
parameters $t_7$ and $t_8$. After some tedious computations, we obtained
\begin{eqnarray}\label{key-eheq}
&&
\frac{\mathcal{D}(\underline{t}') U(\underline{t})}{\theta\left(\frac{t_3}{t_3'},
\frac{t_3}{t_4'};p\right)}
=t_4t_6\frac{\theta\left(\frac{t_3t_4t_7t_8}{q^2},t_7t_6,t_8t_6,\frac{t_7t_8}{q};p
\right)}
{\theta\left(\frac{t_3t_7}{q},\frac{t_3t_8}{q},\frac{t_4t_7}{q},\frac{t_4t_8}{q},
\frac{qt_6}{t_5};p\right)}
\\ &&\makebox[2em]{} \times
\prod_{j=1,\, \neq 5,6}^8\theta\left(\frac{t_5t_j}{q};p\right)
U(q^{-1}t_5,qt_6)
-\frac{\theta\left(\frac{qt_6}{t_8};p\right)\prod_{j=1}^4\theta(t_6t_j;p)}
{\theta\left(\frac{t_3t_8}{q},\frac{q}{t_4t_8};p\right)}
\nonumber \\ &&\makebox[2em]{} \times
\left( \frac{\theta\left(t_6t_8,\frac{t_7}{t_5},\frac{t_5t_7}{q},
\frac{t_3t_4t_7t_8}{q^2};p\right)}
{\theta\left(\frac{qt_6}{t_5},\frac{t_3t_7}{q},
\frac{t_4t_7}{q};p\right)} -
\frac{\theta\left(\frac{t_7t_8}{q},t_5t_6,\frac{t_3t_4t_6t_8}{q};p\right)}
{\theta\left(t_3t_6,t_4t_6;p\right)}  \right)
 U(\underline{t}).
\nonumber\end{eqnarray}
The elliptic hypergeometric equation easily follows from this relation.
Indeed, it is sufficient to replace $t_3'$ by $t_3''$ in  \eqref{key-eheq},
to subtract resulting equation from \eqref{key-eheq},
and to use operator identity \eqref{op-ident}.

In the new parameter notation
$$
\ve_k=\frac{q}{ct_k},\; k=1,\ldots,5,\quad
\ve_6= c t_6 p^4, \quad \ve_8=\frac{c}{t_6},
\quad \ve_7=\frac{\ve_8}{q}
$$
with the balancing condition $\prod_{k=1}^8\ve_k=p^2q^2$,
the operator $\mathcal{D}$ takes the form
\be
\mathcal{D}_x(\ve_1,\ldots,\ve_8):= A(x)(T_{x,q}-1)+ A(x^{-1})(T_{x,q}^{-1}-1)+ \nu,
\label{D}\ee
where the coefficients
\be
A(x)=\frac{\prod_{k=1}^8 \theta(\ve_kx;p)}{\theta(x^2,qx^2;p)},
\qquad
\nu=\prod_{k=1}^6\theta\left(\frac{\ve_k\ve_8}{q};p\right)
\lab{pot}\ee
obey the explicit $S_6$-symmetry in $\ve_1,\ldots,\ve_6$.
Evidently, the equation $\mathcal{D}_x(\underline{\ve})\psi(x;\underline{\ve})=0$
has a solution $\psi(x;\underline{\ve})=U(t_1,\ldots,t_6,cx,cx^{-1})$, where
$$
t_k=\frac{q}{c\ve_k},\; k=1,\ldots,5,\quad
t_6=\frac{c}{\ve_8},\quad c=\frac{\sqrt{\ve_6\ve_8}}{p^2}.
$$
The operator $\mathcal{D}_x(\underline{\ve})$
has the following conjugate with respect to the inner product
$\langle \psi,\phi\rangle$ \eqref{sp}:
\begin{eqnarray}
\mathcal{D}_x^*(\ve_1,\ldots,\ve_8)=\alpha(x)\mathcal{D}_x(\ve_1,\ldots,\ve_8)\alpha(x)^{-1},
\quad
\alpha(x)=\prod_{k=1}^8\Gamma_{\! p,q}(\ve_k x^{\pm}),
\label{conj}\end{eqnarray}
provided the functions $\psi$ and $\phi$ do not have singularities
in the annulus $|q|<|x|<|q|^{-1}$.
A solution of the equation $\mathcal{D}_x^*(\ve)\chi(x;\ve)=0$ is therefore
given by the function
$$
\chi(x;\ve)=\alpha(x)U(t_1,\ldots,t_6,cx^{-1},cx).
$$
We denote as $\underline{\ve}'$ the set of parameters $\underline{\ve}$
where $\ve_{3,4}$ are replaced with $\ve_{3,4}'$ such that
$\ve_3\ve_4=\ve_3'\ve_4'$ (or, equivalently, $t_{3,4}\to t_{3,4}'$,
$t_3t_4=t_3't_4'$). The evident equality
\begin{equation}
0=\langle \mathcal{D}_x^*(\underline{\ve}')
\chi(x;\underline{\ve}'),\beta(x)
\psi(x;\underline{\ve})\rangle=\langle \chi(x;\underline{\ve}'),
\mathcal{D}_x(\underline{\ve}')\beta(x)\psi(x;\underline{\ve})\rangle,
\label{eheq-bio}\end{equation}
where $\beta(x)$ is some function, generates an infinite family
of biorthogonality relations. For instance, if $\beta=1$,
then the function $\chi(x;\underline{\ve}')$ is orthogonal
to the function on the right-hand side of \eqref{key-eheq} for
$\theta(\ve_3'/\ve_3,\ve_4'/\ve_3)\neq0$. This complicated
biorthogonality is clearly different from \eqref{V-ort}.
However, it is plausible that for some particular
$\beta(x)$ relation \eqref{eheq-bio} yields equality \eqref{V-ort}.

\section{Connection to the Sklyanin algebra}

The Sklyanin algebra \cite{skl1,skl2} is generated by four operators
$S_0,\ldots,S_3$ satisfying relations
\begin{eqnarray}\nonumber
&& S_\alpha S_\beta -S_\beta S_\alpha =i(S_0S_\gamma+S_\gamma S_0),
\\ &&
S_0S_\alpha - S_\alpha S_0 = iJ_{\beta\gamma}(S_\beta S_\gamma+S_\gamma S_\beta),
\label{s-rel2}\end{eqnarray}
where the triple $(\alpha,\beta,\gamma)$ represents any cyclic
permutation of $(1,2,3)$.

Sklyanin has found the following realization of this algebra in
the space of theta functions of the even order $2N$:
$$
S_a=\frac{1}{\theta_1(2u)}\left(s_a(u-g)e^{\eta \partial_u}
-s_a(-u-g)e^{-\eta \partial_u}\right), \quad g:=\frac{N\eta}{2},
$$
where $e^{\pm\eta \partial_u}$ are the shift operators,
$e^{\pm\eta \partial_u}f(u)=f(u\pm\eta)$, and
\begin{eqnarray*}
&& s_0(u)=\theta_1(\eta,2u), \quad
s_1(u)=\theta_2(\eta,2u)=\theta_1(\eta+1/2,2u+1/2),
\\ &&
s_2(u)=i\theta_3(\eta,2u)=ie^{\pi i(\tau/2+\eta+2u)}
\theta_1(\eta+1/2+\tau/2,2u+1/2+\tau/2),
\\ &&
s_3(u)=\theta_4(\eta,2u)=-e^{\pi i(\tau/2+\eta+2u)}
\theta_1(\eta+\tau/2,2u+\tau/2).
\end{eqnarray*}
The structure constants $J_{\alpha\beta}$ in this case are parametrized as
$$
J_{12}=\frac{\theta_1^2(\eta)\theta_4^2(\eta)}
{\theta_2^2(\eta)\theta_3^2(\eta)}, \quad
J_{23}=\frac{\theta_1^2(\eta)\theta_2^2(\eta)}
{\theta_3^2(\eta)\theta_4^2(\eta)}, \quad
J_{31}=-\frac{\theta_1^2(\eta)\theta_3^2(\eta)}
{\theta_2^2(\eta)\theta_4^2(\eta)}.
$$
They can be represented in the form $J_{\alpha\beta}=(J_\beta-J_\alpha)/J_\gamma$,
where
$$
J_1=\frac{\theta_2(2\eta)\theta_2(0)}{\theta_2^2(\eta)},\quad
J_2=\frac{\theta_3(2\eta)\theta_3(0)}{\theta_3^2(\eta)},\quad
J_3=\frac{\theta_4(2\eta)\theta_4(0)}{\theta_4^2(\eta)}.
$$

Verification of the defining relations is rather complicated and
requires multiple use of the Jacobi identity \cite{ros:sklyanin}
\begin{equation}
2\theta_1(\vec b-B/2)=\theta_1(\vec b)+\theta_1(\vec b+1/2)+e^{\pi i(\tau+B)}
(\theta_1(\vec b+\tau/2)-\theta_1(\vec b+1/2+\tau/2)),
\label{jac-id}\end{equation}
where $B=b_1+b_2+b_3+b_4$ and $\theta_1(\vec b)=\prod_{k=1}^4\theta_1(b_k)$.

There are two Casimir operators for this algebra.
The first one is
\begin{eqnarray*}
&& K_0=\sum_{a=0}^3 S_a^2
=\sum_{a=0}^3 (-1)^{\delta_{a,2}} \frac{\theta_{a+1}^2(\eta)}{\theta_1(2u,2u+2\eta)}
\Bigl( \theta_{a+1}(2u-2g,2u-2g+2\eta)e^{2\eta\partial_u}
\\ && \qquad
-\theta_{a+1}(2u-2g,-2u-2g-2\eta) \Bigr)+(u\to -u)
= 4\theta_1^2(2g+\eta).
\end{eqnarray*}
The terms proportional to $e^{\pm 2\eta\partial_u}$ cancel due to the Jacobi
identity with the choice $b_1=b_2=\eta,b_3=2u-2g,b_4=2u-2g+2\eta$,
so that $B=4(u-g+\eta)$ and $\theta_1(\vec b-B/2)=0$. The scalar terms
combine to $4\theta_1^2(2g+\eta)$ due to the same identity with the
choice $b_1=b_2=\eta,b_3=2u-2g,b_4=-2u-2g-2\eta$ with $B=-4g$.

The second Casimir operator has the form
\begin{eqnarray*}
&& K_2=\sum_{a=1}^3 J_aS_a^2
=\sum_{a=1}^3 (-1)^{\delta_{a,2}} \frac{\theta_{a+1}^2(2\eta,0)}
{\theta_1(2u,2u+2\eta)}
\Bigl( \theta_{a+1}(2u-2g,2u-2g+2\eta)e^{2\eta\partial_u}
\\ && \qquad
-\theta_{a+1}(2u-2g,-2u-2g-2\eta) \Bigr)+(u\to -u)
= 4\theta_1(2g,2g+2\eta).
\end{eqnarray*}
The terms proportional to $e^{\pm 2\eta\partial_u}$ cancel due to the Jacobi
identity with the choice $b_1=0, b_2=2\eta,b_3=2u-2g,b_4=-2u+2g-2\eta$,
so that $B=0$ and $\theta_1(\vec b-B/2)=0$. The scalar terms
combine to $4\theta_1(2g,2g+2\eta)$ due to the same identity with the
choice $b_1=0,b_2=2\eta,b_3=2u-2g,b_4=-2u-2g-2\eta$ with $B=-4g$.
We stress that all the described operator relations are satisfied
for arbitrary continuous parameter $g$ (i.e., the choice
$g=N\eta/2$ of \cite{skl2} corresponds to a special finite
dimensional representation in the space of theta functions).

As noticed by Rains \cite{rai:abelian,ros:sklyanin}, a linear combination
of the Sklyanin algebra generators $S_a$ is equivalent to the $D$-operator.
Indeed, it is related to the combination
\begin{eqnarray*}
&& 2\Delta(a_1,a_2,a_3,a_4):=\frac{\prod_{j=1}^3\theta_1(a_j+a_4+2g)}
{\theta_1(\eta)}S_0
-\frac{\prod_{j=1}^3\theta_1(a_j+a_4+2g+\frac{1}{2})}
{\theta_1(\eta+\frac{1}{2})}S_1
\\ && \qquad
-ie^{\pi i(\frac{\tau}{2}+2a_4+2g-\eta)}
\frac{\prod_{j=1}^3\theta_1(a_j+a_4+2g+\frac{1+\tau}{2})}
{\theta_1(\eta+\frac{1+\tau}{2})}S_2
\\ && \qquad
+e^{\pi i(\frac{\tau}{2}+2a_4+2g-\eta)}
\frac{\prod_{j=1}^3\theta_1(a_j+a_4+2g+\frac{\tau}{2})}{\theta_1(\eta+\frac{\tau}{2})}
S_3,
\end{eqnarray*}
where we assume that $a_1+a_2+a_3+a_4=-4g$. Substituting explicit expressions
for  $S_a$ and using the Jacobi identity
with $b_1=a_1+a_4+2g, b_2=a_2+a_4+2g, b_3=a_3+a_4+2g, b_4=-2(u-g)$
(so that $B=2a_4+4g-2u$), we find
$$
\Delta(a_1,a_2,a_3,a_4)=\frac{\prod_{j=1}^4\theta_1(a_j+u)}
{\theta_1(2u)}e^{\eta\partial_u}+\frac{\prod_{j=1}^4\theta_1(a_j-u)}
{\theta_1(-2u)}e^{-\eta\partial_u}.
$$
Conversely, any $S_a$ can be viewed as a special case of the $\Delta$-operator:
\begin{eqnarray*}
&& S_0=\chi \theta_1(\eta)
\Delta\left(-g,\frac{1}{2}-g,\frac{\tau}{2}-g,-\frac{1+\tau}{2}-g\right), \quad
\\ &&
S_1=-\chi \theta_1\left(\eta+\frac{1}{2}\right)\Delta\left(\frac{1}{4}-g,-\frac{1}{4}-g,
\frac{1}{4}+\frac{\tau}{2}-g,-\frac{1}{4}-\frac{\tau}{2}-g\right),
\\ &&
S_2=\chi e^{\pi i\eta}\theta_1\left(\eta+\frac{1+\tau}{2}\right)
\Delta\left(\frac{1+\tau}{4}-g,\frac{1-\tau}{4}-g,
\frac{\tau-1}{4}-g,-\frac{1+\tau}{4}-g\right), \quad
\\ &&
 S_3=\chi e^{\pi i\eta} \theta_1\left(\eta+\frac{\tau}{2}\right)
\Delta\left(\frac{\tau}{4}-g,-\frac{\tau}{4}-g,
\frac{1}{2}+\frac{\tau}{4}-g,-\frac{1}{2}-\frac{\tau}{4}-g\right),
\end{eqnarray*}
where $\chi=ip^{1/8}/(p;p)_\infty$, $p=e^{2\pi i\tau}$.

In the multiplicative system of notation
$$
(a,b,c,d):=e^{2\pi ia_{1,2,3,4}}, \quad \rho:=abcd= e^{-8\pi i g},\quad
z:=e^{2\pi i u},\quad q:=e^{4\pi i\eta},
$$
the operator  $\Delta$ can be rewritten as
$$
\Delta(a_1,a_2,a_3,a_4)=\left(ip^{1/8}(p;p)_\infty\right)^3e^{4\pi i g}
D(a,b,c,d;p;q),
$$
where the operator $D$ coincides with \eqref{D-oper}.
In the paper \cite{rai:abelian} it was indicated that equality \eqref{ccr}
comprises all commutation relations of the Sklyanin algebra.
The solution of our generalized eigenvalue problem \eqref{gevp2},
fixed in \eqref{basis},
defines thus a special infinite dimensional module for the Sklyanin
algebra with continuous values of the Casimir operators.

Relation \eqref{gevp2} describes the action of the $D$-operator on this
module $f(z;w;a,b;\rho)$ with a special choice of parameters.
To know the action of
the Sklyanin algebra generators in general, we need to consider
\begin{eqnarray}\label{gen-act}
&& D(a,b,c,d;q;p)f(z;w;q^{1/2}e,q^{1/2}h;\rho)
\\ && \makebox[2em]{}
=\frac{z}{\theta(z^2;p)}
\left(\frac{\theta(az,\ldots,dz,\sqrt{\frac{eh}{\rho}}w^\pm z;p)}
{z^2\theta(ez,hz;p)}- (z\to z^{-1})
\right)f(z;w;e,h;\rho),
\nonumber\end{eqnarray}
where $\rho=abcd$.
We apply now a known theta functions identity (see Lemma 6.4 in \cite{ros:elementary})
\begin{eqnarray*}
&& z^{-n-1}\prod_{j=1}^n\theta(a_jz;p)\prod_{k=1}^{n+2}\theta(b_jz;p)
-(z\to z^{-1})
\\ && \makebox[2em]{}
=\frac{(-1)^nz\theta(z^{-2};p)}{a_1\cdots a_n}\sum_{j=1}^n
\frac{\prod_{k=1}^{n+2}\theta(a_jb_k;p)\prod_{l=1,\neq j}^n\theta(a_jz^\pm;p)}
{\prod_{l=1,\neq j}^n\theta(a_j/a_l;p)},
\end{eqnarray*}
where $\prod_{j=1}^na_j\prod_{k=1}^{n+2}b_k=1$. Namely, we choose $n=3$,
$b_1=a,\ldots,b_4=d,b_5=e^{-1}, a_1=h^{-1}, a_{2,3}=w^{\pm1} \sqrt{eh/\rho}$
and find after some simplifications on the right hand side of \eqref{gen-act}
\begin{eqnarray*}
&& \rho h\Biggl( \frac{\theta(\sqrt{\frac{eh}{\rho}}aw,\ldots,
\sqrt{\frac{eh}{\rho}}dw,\sqrt{\frac{eh}{\rho}}e^{-1}w;p)}
{h^2\theta(h\sqrt{\frac{eh}{\rho}} w,w^2;p)}f(z;q^{-1/2}w;qe,h;\rho)
+(w\to w^{-1})
\\ && \makebox[2em]{}
+\frac{\theta(h^{-1}a,\ldots,h^{-1}d,h^{-1}e^{-1};p)}
{\theta(h^{-1}\sqrt{\frac{\rho}{eh}} w^\pm;p)}f(z;w;qe,qh;\rho)\Biggr).
\end{eqnarray*}
Fixing parameters $a_{1,2,3}$ in a different admissible way yields
hidden permutational symmetries of this expression. In particular,
it is symmetric in $e$ and $h$.

The $D$-operator maps thus basis vectors to their linear combinations
with shifted values of parameters $e,h$, and $w$. Using relation \eqref{key-rel-fin}
with an appropriate choice of parameters, we can replace
$f(z;q^{-1/2}w;qe,h;\rho)$ by the action of some integral operator
on the vector $f(z;x;qe,qh;\rho)$ with integration over $x$.
This provides a realization of the $D$-operator as an integral
operator acting on basis vectors with  the parameters $qe$ and $qh$.

In \cite{ros:sklyanin}, Rosengren proved the Sklyanin conjecture about
the reproducing kernel associated with the discrete basis \eqref{dis-bas}
and the Sklyanin invariant measure \cite{skl2}.
We use a different measure \eqref{sp}, i.e. relation \eqref{V-fact}
does not represent a generalization of the similar overlap product for
the elliptic $6j$-symbols in \cite{ros:sklyanin}.

\section{Elliptic analogues of the Faddeev modular double}

Now we would like to discuss the uniqueness of the basis vectors
$f(z;w;a,b;\rho)$. As we saw above, they are fixed up to a constant
multiplier by two generalized eigenvalue problems \eqref{GEVP}
and \eqref{GEVP-perm}. The first equation uses linear combinations
of the standard Sklyanin algebra generators
$$
S_a=i^{\delta_{a,2}}\frac{\theta_{a+1}(\eta|\tau)}{\theta_1(2u|\tau)}
\left(\theta_{a+1}(2u-2g|\tau)e^{\eta\partial_u}
-\theta_{a+1}(-2u-2g|\tau)e^{-\eta\partial_u}\right).
$$
The second equation requires introduction of another copy of this
algebra, obtained from the first one by permutation of
$q$ and $p$ (or of $\tau$ and $2\eta$):
$$
\tilde S_a=i^{\delta_{a,2}}\frac{\theta_{a+1}(\frac{\tau}{2}|2\eta)}
{\theta_1(2u|2\eta)}\left(\theta_{a+1}(2u-2g|2\eta)e^{\frac{\tau}{2}\partial_u}
-\theta_{a+1}(-2u-2g|2\eta)e^{-\frac{\tau}{2}\partial_u}\right).
$$
For that, we should have $|q|<1$ (i.e., Im$(\eta)>0$).
The structure constants $\tilde J_{\alpha\beta}$ are also obtained from
$J_{\alpha\beta}$ by permutation of $\tau$ and $2\eta$.
Evidently, the  $\tilde S_a$-operators are not well defined in the limit
$\text{Im}(\tau)\to +\infty$, i.e. such a way of fixing representation modules
exists only at the elliptic level.

The two sets of operators $S_a$ and $\tilde S_a$
satisfy the following simple cross-commutation relations
\begin{eqnarray*}
&& S_a \tilde S_b=\tilde S_b S_a, \quad a,b\in\{0,3\} \quad
\text{or} \quad a,b\in\{1,2\},
\\
&& S_a \tilde S_b=-\tilde S_b S_a, \quad a\in\{0,3\},\; b\in\{1,2\} \quad
\text{or} \quad a\in\{1,2\},\;b\in\{0,3\}.
\end{eqnarray*}

Algebra \eqref{s-rel2} has many automorphisms \cite{skl2}.
For real or purely imaginary $\eta$ (i.e., $q^*=q^{\pm 1}$), the transformation
$S_a^*=-S_a$ is an involution for imaginary $\tau$ (i.e., $p^*=p$).
Because the parameter $\eta$ should lie in the upper half plane for the
$\tilde S_a$-operators, such an involution exists for both algebras
only for purely imaginary $\eta$.
There are specific involutions existing only for our direct product
of two Sklyanin algebras. For instance, at $\eta^*=-\tau/2$ (i.e., $q^*=p$)
the conditions Im$(\tau)$, Im$(\eta)>0$ are preserved and we obtain
$S_a^*=-\tilde S_a$ and $J_{\alpha\beta}^*=\tilde J_{\alpha\beta}$.

There is a second possibility for fixing solutions of equation
\eqref{GEVP} uniquely. For its description, we use additive notation and
substitute $\tau=\omega_3/\omega_2,\; 2\eta=\omega_1/\omega_2$,
$$
 z=e^{2\pi i \frac{u}{\omega_2}},\quad
a,b,c,d=e^{2\pi i\frac{\alpha,\beta,\gamma,\delta}{\omega_2}},\quad
w=e^{2\pi i\frac{v}{\omega_2}},\quad
\rho=e^{2\pi i\frac{\sigma}{\omega_2}}
$$
and $f(z;w;a,b;\rho)=:h(u;v;\alpha,\beta;\sigma)$ in \eqref{f-eq}, which yields
\be
\frac{h(u+\frac{\omega_1}{2};v;\alpha+\frac{\omega_1}{2},
\beta+\frac{\omega_1}{2};\sigma)}
{h(u-\frac{\omega_1}{2};v;\alpha+\frac{\omega_1}{2},
\beta+\frac{\omega_1}{2};\sigma)}=z^4
\frac{\theta(az^{-1},bz^{-1},w^\pm z^{-1}\sqrt{cd};p)}
{\theta(az,bz,w^\pm z\sqrt{cd} ;p)}.
\lab{f-eq-add}\ee
A solution of this equation valid for $|p|<1$ and $|q|\leq 1$
has the form
\begin{equation}
h(u;v;\alpha,\beta;\sigma)=\frac{G(\alpha/2+\beta/2 -\sigma/2
\pm v \pm u;\omega)} {G(\alpha\pm u,\beta \pm u;\omega)},
\label{mod-basis}\end{equation}
where $G(u;\omega)$ is the modified elliptic gamma function.
It is defined up to the multiplication by an
arbitrary function $\varphi(u)$ such that $\varphi(u+\omega_1)=\varphi(u)$.
The gamma function $G(u;\mathbf{\omega})$ is symmetric in $\omega_1$
and $\omega_2$. It remains a well defined meromorphic function
even for $\omega_1/\omega_2>0$ (i.e., $|q|=1$ or $\eta$ is real).
Demanding that function \eqref{mod-basis} solves
simultaneously equation \eqref{f-eq-add} and its partner obtained
by permuting $\omega_1$ and $\omega_2$ for $\omega_1/\omega_2>0$,
we fix $\varphi(u)=constant$.

In terms of the Sklyanin algebras, we have now two copies of it
with the second algebra being obtained from the first one
by permutation of $\omega_1$
with $\omega_2$. Such a transformation is equivalent to changes
$\eta\to 1/4\eta$, $\tau\to\tau/2\eta$,  and $u\to u/2\eta$ in $S_a$, which yields
$$
\tilde S_a=i^{\delta_{a,2}}
\frac{\theta_{a+1}\left(\frac{1}{4\eta}\Big|\frac{\tau}{2\eta}\right)}
{\theta_1\left(\frac{u}{\eta}\Big|\frac{\tau}{2\eta}\right)}
\left(\theta_{a+1}\left(\frac{u-g}{\eta}\Big|\frac{\tau}{2\eta}\right)
e^{\frac{1}{2}\partial_u} -\theta_{a+1}\left(\frac{-u-g}{\eta}\Big|
\frac{\tau}{2\eta}\right) e^{-\frac{1}{2}\partial_u}
\right).
$$
We have now the following cross-commutation relations for $S_a $ and $\tilde S_b$:
\begin{eqnarray*}
&& S_a \tilde S_b=\tilde S_b S_a, \quad a,b\in\{0,1\}\quad
\text{or}\quad a,b\in\{2,3\},
\\
&& S_a \tilde S_b=-\tilde S_b S_a, \quad a\in\{0,1\},\; b\in\{2,3\} \quad
\text{or} \quad a\in\{2,3\},\;b\in\{0,1\}.
\end{eqnarray*}

Consider the degeneration of these Sklyanin algebra generators in the limit
Im$(\omega_3/\omega_1)$, Im$(\omega_3/\omega_2)\to +\infty$
(i.e., $r,p\to 0$), which is permissible now.
It is not difficult to see that for $p\to0$ expressions
$p^{-1/8}S_{0,1},$ $p^{1/8}(S_2+iS_3)$, and $p^{-3/8}(S_2-iS_3)$
have well defined limits to some finite difference operators \cite{gz}.
Simplifying the latter operators by shifting $u\to u+c$
and taking the limit $\epsilon:=e^{-2\pi i c}\to 0$, we obtain
three operators $E,F,$ and $k$:
\begin{eqnarray*}
&& \lim_{\epsilon,p\to 0}\frac{p^{-1/8}iS_0}{q^{1/4}-q^{-1/4}}=k+k^{-1},\qquad
\lim_{\epsilon,p\to 0}\frac{p^{-1/8}S_1}{q^{1/4}+q^{-1/4}}=k-k^{-1},
\quad
\\ &&
k:=e^{-2\pi ig+\eta\partial_u}, \qquad
k^{-1}:=e^{2\pi ig-\eta\partial_u},
\\ &&
\lim_{\epsilon,p\to0}p^{1/8}\frac{S_2+iS_3}{2\epsilon}
=e^{-2\pi i u}(e^{-\eta\partial_u}-e^{\eta\partial_u})
=:(q^{1/2}-q^{-1/2})F,\quad
\\ &&
\lim_{\epsilon,p\to0}\epsilon\frac{S_2-iS_3}{2p^{3/8}}
=e^{2\pi i u}(e^{4\pi ig-\eta\partial_u}
-e^{-4\pi ig+\eta\partial_u})=:(q^{-1/2}-q^{1/2})E,
\end{eqnarray*}
where $q=e^{4\pi i\eta}$.
The operators $E,F,k$ form the quantum algebra $U_q(sl_2)$ with the
defining relations
$$
kE=q^{1/2}E k,\quad kF=q^{-1/2}F k,\quad
EF-FE=\frac{k^2-k^{-2}}{q^{1/2}-q^{-1/2}}.
$$
The $\tilde S_a$-operators degenerate in a similar way to
\begin{eqnarray*}
&& {\tilde k}^{\pm 1}=e^{\mp \frac{\pi ig}{\eta}\pm\frac{1}{2}\partial_u},
\quad
\tilde F=\frac{e^{-\frac{\pi i u}{\eta} }(e^{-\frac{1}{2}\partial_u}
-e^{\frac{1}{2}\partial_u})} {e^{\frac{\pi i}{2\eta} }
- e^{-\frac{\pi i}{2\eta} }},\quad
\\ &&
\tilde E=\frac{e^{\frac{\pi i u}{\eta} }
(e^{\frac{2\pi ig}{\eta}-\frac{1}{2}\partial_u}
-e^{-\frac{2\pi ig}{\eta}+\frac{1}{2}\partial_u})}
{e^{-\frac{\pi i}{2\eta}}-e^{\frac{\pi i}{2\eta} }},
\end{eqnarray*}
forming the $U_{{\tilde q}^{-1}}(sl_2)$-algebra
$$
\tilde k\tilde E={\tilde q}^{-1/2}\tilde E \tilde k,\quad
\tilde k\tilde F={\tilde q}^{1/2}\tilde F \tilde k,\quad
\tilde E\tilde F-\tilde F\tilde E
=\frac{{\tilde k}^2-{\tilde k}^{-2}}{{\tilde q}^{-1/2}-{\tilde q}^{1/2}},
$$
with $\tilde q=e^{-\pi i/\eta}$.
We obtain thus the Faddeev modular double \cite{fad:mod}.

In \cite{fad:mod}, it was demanded that generators of the
copies of quantum algebras commute with each other.
This demand is satisfied by the explicit realization
of the modular double used in \cite{kls:unitary}.
In our case, not all the generators of two algebras commute with each other:
\begin{equation}
k \tilde F=- \tilde F k,\quad
k \tilde E=- \tilde E k,\quad
\tilde k F=- F\tilde  k,\quad
\tilde k E=- E\tilde  k,
\label{anti}\end{equation}
which seems to contradict this demand. However, introducing
the operator $K=k^2$ and rewriting defining relations
of $U_q(sl_2)$-algebra as
$$
KE=qE K,\quad KF=q^{-1}F K,\quad
E F-F E=\frac{K-K^{-1}}{q^{1/2}-q^{-1/2}},
$$
we obtain the desired commutativity of generators
$E,F,K$ with $\tilde E,\tilde F,\tilde K$, as in \cite{kls:unitary}
(note that $q$ and $\tilde q$ of the paper \cite{kls:unitary} equal to our
parameters $q^{1/2}$ and ${\tilde q}^{-1/2}$, respectively). %
At the elliptic level, such a quadratic change is not natural;
one has to work with the operators $S_a$ and $\tilde S_a$
some of which anticommute with each other.

For $|q|=1$, algebra $U_q(sl_2)$ has the involution
$E^*=-E, F^*=-F, K^*=K$. As noticed by Faddeev \cite{fad:mod},
for $q^*=\tilde q$, the modular double has another involution
permuting algebras in the pair: $E^*=-\tilde E, F^*=-\tilde F, K^*=\tilde K$.
For our direct product of the Sklyanin algebras, in the first
case we have $S_a^*=-S_a$ and $\tilde S_a^*=-\tilde S_a$, provided
$J_{\alpha\beta}^*= J_{\alpha\beta}$ and
$\tilde J_{\alpha\beta}^*=\tilde J_{\alpha\beta}$
(or, equivalently, $|q|=1$ and $p^*=p$). An elliptic
analogue of the Faddeev involution
has the form $S_a^*=-\tilde S_a$, and it exists under the conditions
$J_{\alpha\beta}^*=\tilde J_{\alpha\beta}$. These constraints for
the structure constants are satisfied, if $|2\eta|=1$ (i.e.,
$\eta^*=1/4\eta$ or $q^*=\tilde q$) and $\tau^*=-\tau/2\eta$
(i.e., $p^*=r$). Equivalently, $2\eta=e^{i\phi},\,
\phi\in (-\pi,\pi),$ and $\tau=i\rho e^{i\phi/2}$, $\rho>0$.

Summarizing our consideration, we conclude, that the second way of fixing
functional freedom in the definition of basis vectors provides an elliptic
generalization of the Faddeev modular doubling principle.
The first way of fixing the functional space does not
have an evident trigonometric limit, i.e. at the
elliptic level there are two different ``modular doubles".

\section{Bethe eigenfunctions and relation to the Heun equation}

As shown by Sklyanin, all operators $S_a$ preserve the space of
even theta functions of order $2N$, if
$$
g=\frac{N\eta}{2}.
$$
The operator $\Delta(a_1,\ldots,a_4)$ obeys thus the same property
(its multivariate generalization to the root system $BC_n$ proposed
by Rains \cite{rai:trans} obeys similar property as well).
Therefore one can try to diagonalize it in the space of such functions:
\begin{equation}
\Delta(a_1,\ldots,a_4)\psi(u)=E\psi(u),\qquad
\psi(u)=\prod_{m=1}^{N}\theta_1(u\pm u_m),
\label{D-evp}\end{equation}
$\psi(-u)=\psi(u)$. This equation has the following explicit form
\begin{eqnarray*}
&& \frac{\prod_{k=1}^4\theta_1(a_k+u)}{\theta_1(2u)}
\prod_{m=1}^{N}\theta_1(u\pm u_m+\eta)
+\frac{\prod_{k=1}^4\theta_1(a_k-u)}{\theta_1(-2u)}
\prod_{m=1}^{N}\theta_1(u\pm u_m-\eta)
\\ && \qquad
=E\prod_{m=1}^{N}\theta_1(u\pm u_m), \qquad \sum_{k=1}^4a_k=-2N\eta.
\end{eqnarray*}
If we divide this relation by $\psi(u)$, then we obtain on the left-hand side
two elliptic functions of $u$ of the order $2N+4$ and a constant on the
right-hand side. In order to verify this elliptic functions identity
it suffices therefore to check its validity for $2N+5$ different values of $u$
in the fundamental parallelogram  of periods for some particular choice of
the spectral parameter $E$.

First, it is easily verified for $u=0,1/2,\tau/2,(1+\tau)/2$,
when $\theta_1(2u)=0$. Second, the validity of equation \eqref{D-evp}
at $u=\pm u_j$ yields $N$ Bethe ansatz equations for the roots $u_1,\ldots, u_N$:
\begin{equation}
\prod_{k=1}^4\frac{\theta_1(a_k+u_m)}{\theta_1(a_k-u_m)}
=\prod_{n=1}^N\frac{\theta_1(u_m\pm u_n-\eta)}
{\theta_1(u_m\pm u_n+\eta)},\quad m=1,\ldots,N.
\label{bete}\end{equation}
Actually, it suffices to satisfy any $N-1$ equations among \eqref{bete},
since elliptic functions cannot have only one simple pole in the fundamental
region. As a $(2N+5)$-th value of $u$ for which equation \eqref{D-evp} is valid,
we take $u=a_l$ for some $l$ which fixes the parameter $E$ as
$$
E=\frac{\prod_{k=1}^4\theta_1(a_k+a_l)}{\theta_1(2a_l)}
\prod_{n=1}^N\frac{\theta_1(a_l \pm u_n+\eta)}
{\theta_1(a_l\pm u_n)}.
$$
Apparently, this expression should not depend on the choice of $l=1,2,3,4$,
which is rather easy to check for $N=0$.

One can analyze in a similar way a more general Bethe ansatz
$
\psi(u)=e^{\gamma u}\prod_{k=1}^{2N}\theta_1(u-u_k)
$
for some constant $\gamma$ determined from the eigenvalue problem
$\Delta\psi=E\psi$. There are solutions different from \eqref{D-evp}.
For instance, when the operator $\Delta$ is reduced to $S_0$, then
for $N=2n$ the roots $u_k$ can be chosen in such a way that
$\psi(u)=e^{\gamma u}\prod_{k=1}^n\theta_1(2(u-u_k))$, and the
Bethe equations have solutions with a non-trivial parameter $\gamma$ (see
review \cite{rui:rev} and references therein). All these special
eigenfunctions are expected to be related to some finite difference
operators commuting with $\Delta$ and forming some algebraic relations
with it.

We consider now the limit $\eta\to 0$, corresponding to the transition from
finite differences to differentiations. For that we perform a gauge
transformation
\begin{eqnarray*}
&& \tilde\Delta(a_1,\ldots,a_4)=\chi\rho^{-1}(u)\Delta(a_1,\ldots,a_4)\rho(u)
\\ &&
=\chi\frac{\rho(u+\eta)}{\rho(u)}\frac{\theta_1(\vec a+u)}{\theta_1(2u)}
e^{\eta\partial_u}-\chi\frac{\rho(u-\eta)}{\rho(u)}
\frac{\theta_1(\vec a-u)}{\theta_1(2u)}e^{-\eta\partial_u},
\end{eqnarray*}
where $\chi$ is some constant to be specified later on.
Choosing $\rho(u)$ as a solution of the equation
$$
\rho(u+\eta)=\frac{1}{\chi}\frac{\theta_1(2u)}{\theta_1(\vec a+u)}\rho(u),
$$
we find
$$
\tilde\Delta(a_1,\ldots,a_4)
=e^{\eta\partial_u}-\chi^2\frac{\theta_1(u-\vec a,u+\vec a-\eta)}
{\theta_1(2u,2u-2\eta)}e^{-\eta\partial_u}.
$$
The potential
$$
\frac{\theta_1(u-\vec a,u+\vec a-\eta)}
{\theta_1(2u,2u-2\eta)}
$$
is an elliptic function of $u$ of the eighth order with the periods $1$ and $\tau$.
We choose now
$$
a_1=\alpha_0\eta,\quad a_2=\frac{1}{2}+\alpha_1\eta,\quad
a_3=\frac{\tau}{2}+\alpha_2\eta,\quad a_4=-\frac{1+\tau}{2}+\alpha_3\eta,
$$
and consider the limit $\eta\to 0$. Fixing
$$
\chi=\frac{ip^{1/8}e^{\pi i(\alpha_2-\alpha_3)\eta}}{(p;p)_\infty^3},
$$
we find after some computations $\tilde\Delta(a_1,\ldots,a_4)=2-\eta^2L+O(\eta^3)$,
where
\begin{equation}
L=-\frac{d^2}{du^2}+\sum_{i=0}^3\alpha_i(\alpha_i-1){\mathcal P}\left(u+\frac{\omega_i}{2}
\right)
\label{heun-ham}\end{equation}
and $\omega_0=0, \omega_1=1,\omega_2=\tau, \omega_3=\omega_1+\omega_2$,
and the function
$
{\mathcal P}(u)=-d^2\log \theta_1(u)/du^2
$
coincides up to an additive constant with the Weierstrass function with the periods
$1$ and $\tau$.

The eigenvalue problem for the $L$-operator is known to lead to the
Heun equation with the general set of parameters. The equation determining
the spectrum of the $\Delta$-operator represents thus a finite difference
analog of the Heun equation.
We call therefore $\Delta(a_1,\ldots,a_4)$ (or $D(a,b,c,d;p;q)$)
the difference Heun operator.

Since the Heun equation is related to the one particle Hamiltonian of
the Inozemtsev model \cite{ino:lax} (an elliptic Calogero-Sutherland type
model for the $BC_n$-root system). Therefore is it natural to expect
that the $D$-operator has a similar relation to relativistic quantum
multiparticle systems. Indeed, this is true for the model proposed by
van Diejen \cite{die:integrability} and investigated in detail
by Komori and Hikami \cite{kom-hik:quantum}. In the presence of
the balancing condition $\prod_{j=1}^8 \ve_j=p^2q^2$, the corresponding
one particle Hamiltonian takes the form
\begin{equation}
H=A(x)(T_{x,q}-1)+A(x^{-1})(T_{x,q}^{-1}-1), \quad
 A(x)=\frac{\prod_{m=1}^8\theta(\ve_mx;p)}{\theta(x^2,qx^2;p)}.
\label{vd-ham}\end{equation}
Imposing the additional constraints
$$
\ve_1=q^{1/2},\quad \ve_2=-q^{1/2},\quad
\ve_3=(pq)^{1/2},\quad  \ve_4=-(pq)^{1/2}
$$
and replacing $\ve_5$ with $p\ve_5$, we find $\prod_{j=5}^8\ve_j=1$ and
$$
H=-\ve_5^{-1}\left(
D(\ve_5,\ldots,\ve_8;p;q^2)
-\ve_5^{-1}\theta(\ve_5\ve_6,\ve_5\ve_7,\ve_5\ve_8;p)\right).
$$
A more detailed consideration of relations between this model and
the $D$-operator lies beyond the scope of the present paper.
We see thus that the difference Heun operator is a rather universal object:
it emerges in the Calogero-Sutherland type models,
it generates the Sklyanin algebra and a family of biorthogonal
functions related to the elliptic analogue of the Gauss
hypergeometric function.

The standard eigenvalue problem for Hamiltonian \eqref{vd-ham}
with a special choice of parameters and for a special eigenvalue
coincides with the elliptic hypergeometric equation \cite{spi:cs}.
Therefore, a particular eigenfunction of the $D$-operator
is related to the $V$-function. To find it, we substitute
$t_{1,2}=\pm q^{1/2}c^{-1},\, t_{3,4}=\pm (pq)^{1/2}c^{-1}$ in
\eqref{D-eheq} and find that under these conditions
\begin{eqnarray*}
&& a(t_1,\ldots,t_7,t_8)+a(t_1,\ldots,t_8,t_7)=\kappa(t_1,\ldots,t_8),
\\ &&
\mathcal{D}(t_1,\ldots,t_8)=-\frac{t_6}{c}
D\left(ct_6,\frac{c}{t_6},\frac{c}{qt_6},\frac{pq}{ct_5};p;q^2
\right).
\end{eqnarray*}
This means that the function
$$
\psi_1(x)=\frac{V\left(\frac{q^{1/2}}{c},-\frac{q^{1/2}}{c},
\frac{(pq)^{1/2}}{c},-\frac{(pq)^{1/2}}{c},t_5,t_6,cx^{-1},cx\right)}
{\Gamma_{\! p,q}(cx^\pm t_6^\pm)},
$$
where $t_5t_6=pc^2$, is a zero mode of the resulting $D$-operator. However,
this zero mode is easily computable from scratch
$$
\psi_2(x)=\frac{\Gamma_{\! p,q^2}(\frac{q^2t_6}{c}x^\pm,
\frac{ct_5}{p}x^\pm)}{\Gamma_{\! p,q^2}(qct_6x^\pm,\frac{qc}{t_6}x^\pm)}.
$$
Meromorphic solutions of the equation $D\psi(x)=0$ are defined up to a
multiplication by an arbitrary
function satisfying equality $\varphi(q^2x)=\varphi(x)$,
i.e. the ratio $\psi_1(x)/\psi_2(x)$ is an elliptic function of $x$.

The Heun equation is the general
differential equation of the second order with four
regular singular points. Description of its general solution in
the closed form is one of the important open mathematical problems.
From our consideration we see that solutions of the generalized
eigenvalue problem for the difference Heun operator defines the basis
of functions, the connection problem for which is solved by the
$R$-function \eqref{R}.

The modular parameter $\tau$ can be considered as a time variable for
the Hamilton dynamics generated by the classical analogue of
operator \eqref{heun-ham}. Then the equations of motion
are known to be equivalent to the Painlev\'e-VI equation \cite{man:painleve}.
In a similar way, we can take the classical analogue of the $D$-operator
and interpret it as a non-stationary Hamitonian with the time variable $\tau$.
For a special dependence of the parameters $a_j$ on $\tau$,
the corresponding equations of motion would define
some particular $q$-deformation of the Painleve-VI function.

Biorthogonality \eqref{V-delta} is easily generalizable to a
similar relation for the entire hierarchy of integrals \eqref{ell-hyp-int}.
Therefore, it looks possible to generalize the construction
described in this paper to $I^{(m)}$ with $m>1$.
It is natural to conjecture also that there exist analogues of the basis
vectors $f(z;w;a,b;\rho)$ for the root system $BC_n$, which can be
built with the help of the generalized eigenvalue problem for the
Rains operator \cite{rai:trans}.
Their determination is a challenging problem for future investigations.

\smallskip

{\em Acknowledgments.}
The author is indebted to V.B. Priezzhev
for stimulating discussions and to A.N. Kirillov, E.M. Rains,
H. Rosengren for useful remarks to this paper.
Relations \eqref{key-rel-fin} and \eqref{gevp2}
were presented for the first time in a slightly different form
at the conference {\em Riemann-Hilbert Problems, Integrability and
Asymptotics} (September 20-25, 2005, SISSA, Trieste). The organizers
of this conference are thanked for warm hospitality.
This work is supported in part by the
Russian Foundation for Basic Research (RFBR),
grant no. 06-01-00191.

After completion of this work, L. D. Faddeev communicated
to the author that an anticommutativity related to the
$k, \tilde k$-operators \eqref{anti} was also noticed in
the paper L. D. Faddeev, A. Yu. Volkov, {\em Discrete evolution
for the zero-modes of the Quantum Liouville Model}, arXiv:0803.0230.

\end{document}